\documentclass[graybox]{svmult}

\usepackage{helvet}         
\usepackage{courier}        
\usepackage{type1cm}        
\usepackage{makeidx}         
\usepackage{graphicx}        
\usepackage{multicol}        
\usepackage[bottom]{footmisc}

\usepackage{amsfonts}
\usepackage{amssymb}
\usepackage{amsmath}

\usepackage{cite,bm,comment}
\usepackage{epsf}

\def\beq{\begin{equation}}
\def\eqn#1{\beq\label{#1}}
\def\eeq{\end{equation}}
\def\ee{\end{equation}}

\def\bea{\begin{eqnarray}}
\def\eqnn#1{\begin{eqnarray}\label{#1}}
\def\eea{\end{eqnarray}}

\newcommand{\eqna}[1]{\begin{subequations} \label{#1}
\begin{eqnarray}}
\def\eena{\end{eqnarray}
\end{subequations}}

\def\htt{\hat{T}}
\def\hV{\hat{V}}

\def\rank{{\rm rank}}
\def\hc{\hat{C}}
\def\qrt{{\textstyle{1\over4}}}

\def\nd{\end{document}}
\def\rg{\rangle} \def\lg{\langle}
\def\ha{{\textstyle{1\over2}}}
\def\bbr{\mathbb{R}}
\def\hp{\hat{\varphi}}

\def\vf{\varphi}
\def\D{\Delta} 
 \def\l{\lambda} \def\z{\zeta}
  \def\s{{\sigma}}
\def\Ga{\Gamma}

\def\nn{\nonumber}

\def\del#1{ \partial_{#1} }
\def\s#1{{\mathfrak s}_{#1}}
\def\Del#1{ \frac{\partial}{\partial #1} }
 \def\cg{{\cal G}}

\def\np{\newpage}

\def\ra{\rightarrow}

\def\ca{{\cal A}}  \def\cc{{\cal C}}
  
\def\cg{{\cal G}} \def\ch{{\cal H}} 
 \def\ck{{\cal K}} 
\def\cm{{\cal M}} \def\cn{{\cal N}} 
\def\cp{{\cal P}} \def\cq{{\cal Q}}

\def\tcn{\widetilde{{\cal N}}}
\def\tN{\widetilde{N}}
\def\tn{\widetilde{n}}

\def\tl{{\tilde{\l}}}
\def\t{\tau}

\def\tD{{\tilde{D}}}

\begin{document}

\title*{Anti de Sitter  Holography via Sekiguchi Decomposition}

\author{Vladimir K. Dobrev and Patrick Moylan}
\authorrunning{Dobrev and Moylan}
\institute{Vladimir K. Dobrev  \at
 Institute for Nuclear Research and Nuclear Energy, Bulgarian Academy of
 Sciences, Tsarigradsko Chaussee 72, BG-1784 Sofia, Bulgaria,
 \email{dobrev@inrne.bas.bg}, (corresponding author) \and Patrick Moylan \at Department of  Physics, Pennsylvania
State University, The Abington College, Abington, PA 19001, USA,
\email{pjm11@psu.edu}}

 \maketitle

\abstract{In the present paper we start consideration of anti de
Sitter  holography in the general case of the $(q+1)$-dimensional
anti de Sitter bulk with boundary $q$-dimensional Minkowski
space-time. We present the group-theoretic foundations that are
necessary in our approach. Comparing what is done for $q=3$ the new
element in the present paper is the presentation of the bulk space
as the homogeneous space ~$G/H = SO(q,2)/SO(q,1)$, which homogeneous
space was studied   by Sekiguchi.}


\section{Introduction}

For the last fifteen years due to the remarkable proposal of
\cite{Malda} the AdS/CFT correspondence is a dominant subject in
string theory and conformal field theory. Actually the possible
relation of field theory on anti de Sitter space to conformal field
theory on boundary Minkowski space-time was studied also before,
cf., e.g., \cite{FlFr,AFFS,Fro,BrFr,NiSe,FeFr}. The proposal of
\cite{Malda} was further elaborated in \cite{GKP} and \cite{Wita}.
After that  there was an explosion of related research which
continues also currently.

Let us recall that the AdS/CFT correspondence
has 2 ingredients \cite{Malda,GKP,Wita}:
1. the holography principle, which is very
old, and  means the reconstruction of some objects in the bulk (that
may be classical or quantum) from some objects on the boundary; 2.
the reconstruction of quantum objects, like 2-point functions on the
boundary, from appropriate actions on the bulk.

Our focus is on the first ingredient.
We note that until recently the explicit presentation of
the holography principle was realized in the Euclidean case, i.e., for the group ~$SO(q+1,1)$~
relying on Wick rotations of the final results, cf., e.g. \cite{Wita,Dobads}.
Yet it is desirable to  show the holography principle
by direct construction in Minkowski space-time, i.e., for the conformal
group ~$SO(q,2)$.

This was done for the case $q=3$ in detail in \cite{AizDob}.
In the present paper we start consideration of the general
case of the $(q+1)$-dimensional anti de Sitter bulk with
boundary $q$-dimensional Minkowski space-time.
Actually, here we only lay the group-theoretic foundations that are necessary in our approach
while the actual construction  is postponed to \cite{DoMo2}. As historical remark we mention that
this approach originated in the
construction of the discrete series of unitary representations in
\cite{Hot,Schm}, which was then applied in \cite{DMPPT} for the Euclidean conformal group $SO(4,1)$.
A different approach was applied to the general
 Euclidean case $SO(N,1)$ in \cite{Dobads}. Also  the nonrelativistic Schr\"odinger
  algebra case was considered in \cite{AiDo}.

The new element in the present paper is the presentation of the bulk space as the homogeneous space
~$G/H = SO(q,2)/SO(q,1)$.
For this we use the Sekiguchi decomposition \cite{Seki}
$$ G \cong |_{\rm loc}\ \tN A H $$ 
where $A$ is the subgroup of dilatations, $\tN$ is isomorphic to the
subgroup of translations. The above means that the subgroup $\tN A
H$ is an open dense set of $G$, and thus the homogeneous space
~$G/H$~ is locally isomorphic to bulk space ~$\tN A$.


\section{Preliminaries}


We need some well-known preliminaries to set up our notation and
conventions. The Lie algebra ~$\cg ~=~  so(q,2)$~ may be defined as
the set of $(q+2)\times (q+2)$ matrices ~$X$~ which fulfil the
relation:  \beq
  {}^tX\eta + \eta X = 0,
  \label{so32condition}
\eeq
where  the metric $\eta$ is given by
\beq
  \eta = (\eta_{AB}) = \rm{diag}(-1,1,\ldots,1,-1), \quad A, B = 0, 1, \cdots,
  q+1
  \label{so32metric}
\eeq
Then we can choose a basis $ X_{AB} = - X_{BA} $ of $\cg$ satisfying the commutation relations
\beq
  [X_{AB}, X_{CD}] = \eta_{AC} X_{BD} + \eta_{BD} X_{AC} - \eta_{AD} X_{BC} - \eta_{BC} X_{AD}.
  \label{so32com}
\eeq

We list the important subalgebras of $\cg$:
\begin{itemize}
\item $\ck = so(q) \oplus so(2)$, ~generators: ~$X_{AB}$ ~:~ $(A,B) \in  \{1,\ldots,q\}, \{0,q+1\}$, ~~maximal
compact subalgebra;
\item $\cq$, ~generators: ~$X_{AB}$ ~:~ $A \in  \{1,\ldots,q\},\ B\in \{0,q+1\}$, 
non-compact completion of $\ck$;
\item $\ca = so(1,1)$, ~generator: ~$D \doteq X_{q,q+1}\,$, ~dilatations;
\item $\cm = so(q-1,1)$, ~generators: ~$X_{AB}$ ~:~ $(A,B) \in  \{0,\ldots,q-1\}$, ~Lorentz subalgebra;
\item $\cn$, ~generators: ~$T_{\mu} = X_{\mu q} + X_{\mu,q+1 }$, ~$\mu = 0,\ldots,q-1$, ~translations;
\item $\tcn$,  ~generators: ~$C_{\mu} = X_{\mu q} - X_{\mu,q+1}$, ~$\mu = 0,\ldots,q-1$, ~special conformal
transformations;
\item $\ca_0 = so(1,1)\oplus so(1,1)$, ~generators: ~$X_{0,q-1}, X_{q,q+1}\,$;
\item $\cm_0 = so(q-2)$, ~generators: ~$X_{AB}$ ~:~ $(A,B) \in  \{1,\ldots,q-2\}$;
\item $\cn_0$, ~generators: ~$T_{\mu}$, ~$\mu = 0,\ldots,q-1$, ~
~$T'_{\mu} = X_{\mu 0} + X_{\mu,q-1 }$, ~$\mu = 1,\ldots,q-2$,
 ~extended translations;
 \item $\tcn_0$, ~generators: ~$C_{\mu}$, ~$\mu = 0,\ldots,q-1$, ~
~$C'_{\mu} = X_{\mu 0} - X_{\mu,q-1 }$, ~$\mu = 1,\ldots,q-2$,
 ~extended special conformal transformations;
\item $\ch = so(q,1) $, ~generators: ~$X_{AB}$ ~:~ $(A,B) \in
\{0,\ldots,q\}$.
\end{itemize}
The last subalgebra is the analog of the maximal compact subalgebra
~$so(q+1)$~ of the Euclidean conformal algebra ~$so(q+1,1)$~ of
$q$-dimensional Euclidean space. Thus, it may result from the Wick
rotation of the Euclidean conformal algebra ~$so(q+1,1)$~ to the
Minkowskian conformal algebra ~$so(q,2)$.

Thus, we have several decompositions:
\begin{itemize}
\item $\cg = \ck \oplus \cq$, ~Cartan decomposition;
\item $\cg = \ck \oplus \ca_0 \oplus \cn_0$, ~(also ~$\cn_0 \ra\tcn_0$), ~Iwasawa decomposition;
\item $\cg = \cn_0 \oplus \cm_0 \oplus \ca_0 \oplus \tcn_0$, ~minimal Bruhat
decomposition;
\item $\cg = \cn \oplus \cm \oplus \ca \oplus \tcn$, ~maximal Bruhat
decomposition;
\item $\cg = \ch \oplus \ca \oplus \cn$,  ~(also ~$\cn \ra\tcn$), ~Sekiguchi decomposition
\cite{Seki}.
\end{itemize}
The subalgebra ~$\cp_0 = \cm_0 \oplus \ca_0 \oplus \tcn_0$~ is a
minimal parabolic subalgebra of ~$\cg$. The subalgebra ~$\cp = \cm
\oplus \ca \oplus \tcn$~ is a maximal parabolic subalgebra of
~$\cg$.

Finally, we introduce the corresponding   Lie groups:\\
~$G ~=~ SO_0(q,2)$ with Lie algebra ~$\cg=so(q,2)$, ~$H ~=~
SO(q,1)$ with Lie algebra ~$\ch=so(q,1)$, ~$K ~=~ SO(q) \times
SO(2)$~  is the maximal compact subgroup of $G$, ~$A_0 ~=~ \exp
(\ca_0) ~=~ SO_0(1,1)\times SO_0(1,1) $~ is abelian simply connected, ~$N_0 ~=~ \exp
(\cn_0) ~\cong~ \tN_0 ~=~ \exp (\tcn_0)$,~ are abelian simply connected
subgroups of ~$G$~ preserved by the action of ~$A_0$. The group  ~$M_0
~\cong~ SO_0(q-2)$~ (with Lie algebra $\cm_0$) commutes with $A_0$. Further
$A ~=~ \exp (\ca) ~=~ SO_0(1,1)$~ is abelian simply connected, ~$N
~=~ \exp (\cn) ~\cong~ \tN ~=~ \exp (\tcn)$,~ are abelian simply
connected subgroups of ~$G$~ preserved by the action of ~$A$. The
group  ~$M ~\cong~ SO_0(q-1,1)$~ (with Lie algebra $\cm$) commutes
with $A$.

We mention also some group  decompositions:
\eqna{grdeco}
 G &=& K A_0 N_0, ~({\rm also} ~N_0 \ra \tN_0), ~{\rm Iwasawa~ decomposition};\\
  G &\cong& |_{\rm loc}\
 \tN A  M N , ~{\rm maximal ~Bruhat~ decomposition};\\
 G &\cong& |_{\rm loc}\  \tN A H, ~({\rm also} ~\tN \ra N), ~{\rm Sekiguchi~ decomposition}
\eena In (\ref{grdeco}b,c) the groups on the RHS are open dense
subsets of $G$. We should note that in \cite{Seki} was studied the
more general case ~$SO_0(q,r+1)/SO(q,r)$.

The subgroup ~$P_0 = M_0A_0N_0$~ is a  ~{\it minimal parabolic subgroup}~ of
$G$. The subgroup ~$P = MAN$~ is a  ~{\it maximal parabolic subgroup}~ of
$G$. Parabolic subgroups are important because the representations
induced from them generate all admissible irreducible
representations of semisimple groups \cite{Lan,KnZu}.
The group (algebra) ~$SO_0(q,2)$~ ($so(q,2)$) has one more maximal (cuspidal) parabolic subgroup
(subalgebra) which we do not give here for the lack of space, cf., e.g.,
\cite{Dojmp,DoMo} for $q=4$.

\section{Elementary representations}

We use the approach of \cite{Dob} which we adapt in a condensed form
here. We work with  so-called ~{\it elementary representations}
(ERs). They are induced from representations of the parabolic
subgroups. Here we work with the maximal parabolic ~$P ~=~ MAN\,$,
where we use (non-unitary) finite-dimensional representations ~$\l$~
of $M=SO(q-1,1)$ in the space $V_\l\,$, (non-unitary) characters of
$A$ represented by the conformal weight $\D$, and the factor $N$ is
represented trivially. For further use we give explicit
parametrization of $\l$: \eqn{mlamb} \l = (\l_1\,,\ldots,\l_{\hat
q}) \ , \quad \hat{q} \doteq [\ha q] \ee where $[w]$ is the largest
integer not greater than $w$. The numbers ~$\l_i$~ are all integer or
all half-integer and they fulfill the following conditions:
\eqnn{mlamc} && |\l_1| \leq \l_2 \leq \cdots \leq \l_{q/2} \ , \quad {\rm for} ~ q ~{\rm even} \\
&& 0\leq \l_1 \leq \l_2 \leq \cdots \leq \l_{(q-1)/2} \ , \quad {\rm for} ~ q ~{\rm odd } \nn\eea
 The data ~$\l,\D$~ is enough to determine a weight
~$\chi\in\ch_\cg^*$, where $\ch_\cg$ is the Cartan subalgebra of
$\cg$, cf. \cite{Dob}. Thus, we shall denote the ERs by ~$C^\chi$.
Sometimes we shall write: ~$\chi = [\l,\D]$. The representation
spaces are $C^\infty$ functions on $G/P$, or equivalently, on the
locally isomorphic group ~$N$~ with appropriate asymptotic
conditions (which we do not need explicitly, cf., e.g.
\cite{Dob,DoMo}). We recall that $\tN$ is isomorphic to
q-dimensional Minkowski space-time ~$\mathfrak{M}$~ whose elements
will be denoted by ~$x = (x_0,\ldots,x_{q-1})$, while the
corresponding elements of $\tN$ will be denoted by ~$n_x\,$. The
Lorentzian inner product in ~$\mathfrak{M}$~ is defined as usual:
\beq
  \lg x,x'\rg  \doteq  x_0 x'_0 - \cdots - x_{q-1} x'_{q-1}\ ,
  \label{innerpro}
\eeq
and we use the notation $ \bm{x}^2 = \lg x,x \rg\ . $

The representation action is given as follows:
\eqn{rac} (T^\chi (g)\vf) (x) ~=~ y^{-\D}\, D^\l (m)\, \vf (x') \eeq
the various factors being defined from the local Bruhat decomposition
(\ref{grdeco}b) ~$G \cong_{\rm loc} \tN AMN$~:
\eqn{vars} g^{-1}\,\tn_x ~=~ \tn_{x'}\, a_y^{-1} m^{-1} n^{-1} \ , \eeq
where $\,y\in\bbr_+$ parametrizes the elements ~$a\in A$,
~$m\in M$,  ~$D^\l (m)$~ denotes the representation action of $M$ on the space ~$V_\l$,
~$n\in N$.

On these functions the infinitesimal action of our representations looks as follows:
\bea\label{Boundary}
  & & T_{\mu} = \del{\mu}, \quad \del{\mu} \doteq  \frac{\partial}{\partial x_{\mu}}, \quad \mu = 0, \ldots, q-1,
   \\
  & & D =  -\sum_{\mu=0}^{q-1}  x_{\mu}  \del{\mu} -\Delta,
  \nn \\
   & & X_{0a} = x_0 \del{a} + x_a \del{0} + \s{0a}, \quad a=1,\ldots,q-1,
  \nn \\
  & & X_{ab} = -x_a \del{b} + x_b \del{a} + \s{ab}, \quad 1\leq a<b\leq q-1,
  \nn \\
 & & C_\mu = - 2 \eta_{\mu\mu}\, x_\mu D + \bm{x}^2 \del{\mu} - 2 \sum_{\nu=0}^{q-1} x^\nu \s{\mu\nu} ,
  \nn
 \eea
where ~$\s{\mu\nu}$~ are the infinitesimal generators of ~$D^\l (m)\,$.

We recall several facts about elementary representations \cite{DMPPT,Dob}:
\begin{itemize}
\item The Casimir operators ~$\cc_i$ of $\cg$  have constant values on the ERs:
\eqn{cas} \cc_i(\{ X\})\, \vf(x) ~=~ \chi_i(\l,\D)\, \vf(x) ~, \qquad
i=1,\dots,\rank \,G ~=~ [\ha q]+1, \ee
where ~$\{X\}$~ denotes symbolically the generators
of the Lie algebra ~$\cg$~ of ~$G$,  the action of which
is given in \eqref{Boundary}.
\item
On the ERs are defined the integral Knapp-Stein ~$G_\chi$~
operators which intertwine the representation ~$\chi$~ with the representation
~$\tilde{\chi} \doteq [\tilde{\l}, q-\D]$, where ~$\tilde{\l}$ is the mirror image
of ~$\l$. We recall that the mirror image $\tilde{\l}$ is equivalent to $\l$ when $q$ is odd, while
 for $q$ even and $\l$ parametrized as in \eqref{mlamb}:
~$\l = (\l_1,\l_2,\ldots,\l_{q/2})$~ we have
~$\tilde{\l} = (-\l_1,\l_2,\ldots,\l_{q/2})$.
\item  The representations ~$\chi$~ and
~$\tilde{\chi}$ are called ~{\it partially equivalent}~   due to the existence of the
  intertwining operator $G_\chi$ between them. The representations
are called ~{\it equivalent}~ if the intertwining operator  $G_\chi$
is onto and invertible.
\item
We also recall that   the Casimirs
 $\chi_i$ have the same values on the partially   equivalent ERs:
\eqn{casiv} \chi_i(\l,\D) ~=~ \chi_i(\tl,q-\D) \ee
\end{itemize}

In the above general definition ~$\vf (x)$~ are considered as elements of the
finite-dimensional representation space $V^\l$ in which act the operators $D^\l (m)$.
The representation space $C^\chi$ can be
thought of as the space of smooth sections of the homogeneous vector
bundle (called also vector ~$G$-bundle) with base space ~$G/P$~ and
fibre ~$V_\l\,$, (which is an associated bundle to the principal
~$P$-bundle with total space ~$G$). Actually, we do not need this
description, but following \cite{Dob}  we replace the above
homogeneous vector bundle  with a line bundle again with base space ~$G/P$. The
resulting functions ~$\hp$~ can be thought of as smooth sections of this line
bundle.

In the case when the representation $\l$ is of symmetric traceless tensors of rank $\ell$,
i.e., ~$\l = (0,\ldots,0,\ell)$, we can be more explicit
following \cite{DMPPT}. Namely,
the functions $\hp$ are  scalar functions over an extended space
~$\mathfrak{M} \times ~\mathfrak{M}_0$, where ~$\mathfrak{M}_0$~
is a cone parametrized by the  variable $\z~=~ (\z_0, \ldots, \z_{q-1}) $ subject
to the condition:\beq
  \bm{\z}^2 = \lg \z,\z\rg = \z_0^2 - \cdots - \z_{q-1}^2 = 0.
  \label{conditionforz}
\eeq

The functions on the extended space will be denoted as ~$\hp(x,\z)$.
The   internal variable $\z$ will  carry the representation ~$D^\l$.
Thus, on the functions $\hp$ the infinitesimal generators ~$\s{\mu\nu}$~  from \eqref{Boundary}
are given as follows:
   \beq
  \s{0a} = \z_0 \Del{\z_a} + \z_a \Del{\z_0}, \qquad
      \s{ab} = - \z_a \Del{\z_b} + \z_b \Del{\z_a},
  \label{so12}
\eeq

\section{Bulk representations}

It is well known that the group ~$SO(q,2)$~ is called also anti de
Sitter group, as it is the group of isometry of (q+1)-dimensional
anti de Sitter space: \eqn{adss} \xi^A \,\xi^B\, \eta_{AB} ~=~ 1 \
,\quad A,B =0,\ldots,q+1\ . \eeq There are several ways to
parametrize anti de Sitter space. For $q=3$ in the paper
\cite{AizDob} was  utilized the same local Bruhat decomposition
(\ref{grdeco}b) that we used in the previous section. In the present
paper we shall use the Sekiguchi decomposition (\ref{grdeco}c),
i.e., the factor-space ~$G/H \cong \tN A$. In fact, we use
isomorphic (w.r.t. \cite{AizDob})
 coordinates ~$(x,y) ~=~
(x_0,\ldots,x_{q-1},y)$, ~$y\in\bbr_+\,$. In this setting anti de Sitter
space is called ~{\it bulk}~ space,
while $q$-dimensional Minkowski space-time is called ~{\it boundary}~
space, as it is identified with the
bulk boundary value ~$y=0$.

It is natural to discuss
representations on anti de Sitter space ~$\tN A$~
which are induced from the  subgroup ~$H=SO(q,1)$. Namely, we consider
the representation space:
\eqn{funk} \hc_\t   ~=~ \{ \phi \in
C^\infty(\bbr^q\times\bbr_{>0}\,,\hV_\t) \} \ee
where$\,$ $\t\,$ is an arbitrary finite-dimensional irrep of $H$,$\,$
$\hV_\t\,$ is the finite-dimensional representation space of $\t$,
with representation action:
\eqn{lartu}
(\htt^\t(g)\phi) (x,y) ~=~ \tD^\t (h)\, \phi ( {x'},y') \ee
where the Sekiguchi decomposition is used:
\eqn{nah}  g^{-1} \tn_x a_y ~=~ \tn_{x'} a_{y'} h^{-1}  ~, \quad
g\in G ,\, h\in H, \ n_x, n_{x'} \in N , \, a_y, a_{y'} \in A \ee
and $\tD^\t(k)$ is the representation matrix of $\t$ in $\hV_\t\,$.
For later use we give the parametrization of the relevant subgroups:
\eqn{sekig} H ~=~ \left\{ h=\left[\begin{matrix}h' & 0 \cr 0  & \pm 1  \cr \end{matrix}\right]
\vert h \in SO_0(q,2), \quad h' \in SO_0(q,1), \right\} ~\cong~ SO(q,1) \ee
\eqn{dila} A ~=~ \left\{ a_y =
\left( \begin{matrix}
\mathbb{I}_{q} &0 &0 \cr
0 & \cosh (s)  & \sinh (s)  \cr 0& \sinh (s)   & \cosh (s) \cr \end{matrix}\right)
\vert ~y =e^s, ~~~s\in \mathbb{R}~~\right\}\ee
\eqn{trans}
{\tilde N} = \left\{ {\tilde n}_x  = \left[
\begin{matrix}
\mathbb{I}_{1}   && 0                                     && t                                 &&   t \cr
0                && \mathbb{I}_{q-1}                      &&s^\dagger                         &&  s^\dagger \cr
t               && -s                                     && 1 + \frac{t^2}{2} - \frac{s^2}{2}  &&    \frac{t^2}{2} - \frac{s^2}{2} \cr
-t                 && s                                    && \frac{s^2}{2} - \frac{t^2}{2}    &&  1 +\frac{s^2}{2} - \frac{t^2}{2} \cr
\end{matrix} \right]
\mid (t,s) = \frac{x}{\sqrt{2}} \in \mathbb{R}^q
\right\}
\ee

The infinitesimal generators of \eqref{lartu} are given as follows:
\eqnn{Bulk}
  & & {\hat T}_\mu = \del{\mu}, \quad \mu = 0, \ldots, q-1
   \\
  & & {\hat D}  = -\sum_{\mu=0}^{q-1} x_{\mu} \del{\mu} -y \del{y},
  \nn \\
 & & {\hat X}_{0a} = x_0 \del{a} + x_a \del{0} + \s{0a}, \quad a=1,\ldots,q-1
  \nn \\
  & & {\hat X}_{ab} = -x_a \del{b} + x_b \del{a} + \s{ab}, \quad 1\leq a<b\leq q-1
   \nn\\
& & {\hat C}_\mu =  -2 \eta_{\mu\mu}\,x_\mu D + (\bm{x}^2+y^2) \del{\mu} -
2 \sum_{\nu=0}^{q-1} x^\nu \s{\mu\nu} -2 y \Ga_\mu ,
  \nn
\eea
where ~$\s{\mu\nu},\Ga_\mu$~ are infinitesimal generators of ~$\tD^\t (h)$, such that
(due to the compatibility of $\l$ and $\t$)
~$\s{\mu\nu} = \qrt [\Ga_\mu,\Ga_\nu]$, ~$[\s{\mu\nu},\Ga_\rho]= \eta_{\nu\rho}\Ga_\mu -
\eta_{\mu\rho}\Ga_\nu$.

Note that the realization of $ so(q,2) $ on the boundary given in \eqref{Boundary}
may be obtained from \eqref{Bulk} by replacing ~$y \del{y} \to \Delta$~ and then
taking the limit  $ y \to 0$.

What is important is that, unlike the ERs, the  representations
\eqref{lartu} are highly reducible. Our aim is
  to extract from ~$\hc_\t$~  representations that may be equivalent to ~$C_\chi\,$, $\chi=[\l,\D]$.
The first condition for this is that the ~$M$-representation  ~$\l$~ is
contained in the restriction of the ~$H$-representation ~$\t$~ to ~$M$, i.e.,
~$\l\in\t\vert_M\,$.
Another condition  is that the two representations
would have the same Casimir values ~$\l_i(\l,\D)$.

This procedure is actually well understood and used in the
construction of the discrete series of unitary representations,
cf. \cite{Hot,Schm}, (also \cite{DMPPT,AizDob}$\,$ for $q=4$).
The method  utilizes the fact that in the bulk the Casimir operators
are not fixed numerically. Thus, when a vector-field realization of
the anti de Sitter algebra ~$so(q,2)$~ (e.g., \eqref{Bulk}) is substituted in the bulk
Casimirs the latter turn into   differential  operators.  In contrast, the
boundary Casimir operators are fixed by the quantum numbers of the
fields under consideration. Then the bulk/boundary correspondence
forces  eigenvalue equations involving the Casimir differential
operators. Actually the 2nd order Casimir is enough for this purpose.
That corresponding eigenvalue 2nd order  differential equation is used to find the two-point
Green function in the bulk which is then used to construct the
boundary-to-bulk integral intertwining operator. This operator maps a boundary
field   to a bulk field. For our setting this will be given in detail in \cite{DoMo2}.

Having in mind the degeneracy of Casimir values for partially equivalent
representations (e.g., \eqref{casiv})  we add also the appropriate
asymptotic condition. Furthermore, from now on we shall suppose
that ~$\D$~ is real.

Thus, the representation (partially) equivalent to the ER $\chi$ is defined as:
\eqnn{funb}
 \hc^\t_\chi ~&=&~ \{\, \phi \in \hc_\t ~:~~
\cc_i(\{\hat X\})\, \phi (x,y) ~=~ \l_i(\l,\D)\, \phi (x,y) ~,
\quad \forall i ~, \quad
\l \in\t\vert_M ~, \nn\\
&&\phi (x,y) ~\sim~ y^\D\, \varphi(x) ~ {\rm for}~ y\to 0\ \}
\eea
where $\hat X$ denotes the action \eqref{Bulk} of $\cg$ on the bulk fields.

In the case of symmetric traceless tensors of rank $\ell$ for both ~$M$~  and ~$H$
we can extend the functions on the bulk extended also with the cone $\mathfrak{M}_0\,$.
These extended functions will be denoted by ~$\phi (x,y,\z)$.
On these functions we have the infinitesimal action given by \eqref{Bulk} with
~$\s{\mu\nu}$~ are given by \eqref{so12}, while
 ~$\Ga_\mu$~ are certain finite-dimensional matrices  which we
 shall give in \cite{DoMo2}.

\section{Two parametrizations of bulk space}

As we mentioned in \cite{AizDob} we used as parametrization of the bulk space
the coset ~$G/ M N =|_{\rm loc}\, \tN A$. The local coordinates of this coset come
from the Bruhat decomposition:
\eqn{bruh} g ~=~ \{ g_{AB} \} ~=~ \tn_x  a_y m n \ee
which exists for $g\in G$ forming a dense subset of $G$.
 The local coordinates of the above bulk are:
 \eqnn{bruha}
&& y = \ha (g_{qq} + g_{q,q+1}+g_{q+1,q} + g_{q+1,q+1})~~, \\
&& x_\mu  = \frac{ g_{\mu, q}  +g_{\mu,q+1} }{   g_{qq} + g_{q,q+1}+g_{q+1,q} + g_{q+1,q+1} }
\ , ~~\mu=0,\ldots,q-1
\nn\eea

The parametrization used in the present paper
for the bulk space is the coset ~$G/H =|_{\rm loc}\, \tN A$.
Certainly, it is isomorphic to the bulk above, however,
the local coordinates  are different, namely, the latter .
  come from the Sekiguchi decomposition:
\eqn{sekiguc} g ~=~ \{ g_{AB} \} ~=~ \tn_x   a_y  h \ee
Explicitly, they are given as follows:
 \eqnn{sekigus}
&& y =   | g_{q+1,q} + g_{q+1,q+1} |  \\
&& x_\mu =   \frac{g_{\mu,q+1}} { g_{q+1,q} + g_{q+1,q+1}}\ , ~~\mu=0,\ldots,q-1 \nn\eea

Comparing the two parametrizations \eqref{bruha} and \eqref{sekigus} we see that the latter is
simpler and thus easier to implement. Thus, in the follow-up paper \cite{DoMo2} we shall use
the Sekiguchi decomposition \eqref{sekiguc}.

\vspace{10mm}

\begin{acknowledgement}
The first author has received partial support from COST actions MP-1210
and MP-1405, and from Bulgarian NSF Grant DFNI T02/6.
\end{acknowledgement}

\np


\begin{thebibliography}{99}


\bibitem{Malda}J. Maldacena,~  Adv. Theor. Math. Phys. {\bf 2} (1998) 231
(hep-th/971120).

\bibitem{FlFr}M. Flato and C. Fronsdal,
J. Math. Phys. {\bf 22} (1981) 1100.

\bibitem{AFFS}E. Angelopoulos, M. Flato, C. Fronsdal and D. Sternheimer,
Phys. Rev. {\bf D23} (1981) 1278.

\bibitem{Fro} C. Fronsdal, Phys. Rev. {\bf D26} (1982) 1988.

\bibitem{BrFr}
  P.~Breitenlohner and D.Z.~Freedman,
  Phys. Lett.  {\bf B115} (1982) 197. 

\bibitem{NiSe}H. Nicolai and E. Sezgin, Phys. Let. {\bf 143B} (1984) 103.

\bibitem{FeFr} S. Ferrara and C. Fronsdal,~  Class. Quant. Grav. {\bf 15}
(1998) 2153; ~ 


\bibitem{GKP}S.S. Gubser, I.R. Klebanov and A.M. Polyakov,~
Phys. Lett. {\bf B428} (1998) 105, (hep-th/9802109).

\bibitem{Wita}E. Witten,~  Adv. Theor. Math. Phys. {\bf 2} (1998) 253,
(hep-th/9802150).

\bibitem{Dobads} V.K. Dobrev,
Nucl. Phys. {\bf B553} [PM] (1999) 559. 



\bibitem{AizDob} N. Aizawa and V.K. Dobrev,
 Rept. Math. Phys. {\bf 75} (2015) 179-197.

\bibitem{DoMo2} V.K. Dobrev and P. Moylan, in preparation.

\bibitem{Hot} R. Hotta,
J. Math. Soc. Japan, {\bf 23} (1971) 384.  

\bibitem{Schm} W. Schmid, Rice Univ. Studies, {\bf 56} (1970) 99. 

\bibitem{DMPPT} V.K. Dobrev, G. Mack, V.B. Petkova, S.G. Petrova and I.T.
Todorov, {\it Harmonic Analysis on the $n$ - Dimensional
Lorentz Group and Its Applications to Conformal Quantum Field
Theory}, Lecture Notes in Physics, No 63, 280 pages, Springer
Verlag, Berlin-Heidelberg-New York, 1977.

\bibitem{AiDo} N. Aizawa and V.K. Dobrev,  
Nucl. Phys. {\bf B828} [PM] (2010)  581. 



\bibitem{Seki}J. Sekiguchi, Nagoya Math. J.  {\bf 79} (1980) 151-185.

\bibitem{Lan}R.P. Langlands, {\it  On the classification of irreducible
representations of real algebraic groups}, Math. Surveys and
Monographs, Vol.  31 (AMS, 1988), first as IAS Princeton preprint
(1973).


\bibitem{KnZu} A.W. Knapp and G.J. Zuckerman, in: Lecture Notes in Math.
Vol. 587 (Springer, Berlin, 1977)  pp. 138;
~Ann. Math. {\bf 116} (1982) 389.  

\bibitem{Dojmp}V.K. Dobrev,
~J. Math. Phys. {\bf 26} (1985) 235-251.

\bibitem{DoMo} V.K. Dobrev and P. Moylan, 
Fort. d. Physik, {\bf 42} (1994) 339.  

\bibitem{Dob}V.K. Dobrev,
Rept. Math. Phys. {\bf 25} (1988) 159;
~first as ICTP Trieste preprint IC/86/393 (1986).



\end{thebibliography}
\end{document}